\documentclass[preprint,aps,superscriptaddress,floatfix]{revtex4}
\usepackage{graphicx,epsfig,amssymb,amsmath,array}
\usepackage{enumitem}
\usepackage{relsize,placeins,float}
\usepackage[dvipsnames]{xcolor}
\usepackage{blindtext}
\usepackage[english]{babel}

\begin{document}

\title{Nonlinear delayed forcing drives a non-delayed Duffing oscillator}

\author{Mattia Coccolo}
\affiliation{Nonlinear Dynamics, Chaos and Complex Systems Group, Departamento de F\'{i}sica,
Universidad Rey Juan Carlos, Tulip\'{a}n s/n, 28933 M\'{o}stoles, Madrid, Spain}

\author{Miguel A.F. Sanju\'{a}n}
\affiliation{Nonlinear Dynamics, Chaos and Complex Systems Group, Departamento de F\'{i}sica,
Universidad Rey Juan Carlos, Tulip\'{a}n s/n, 28933 M\'{o}stoles, Madrid, Spain}
\affiliation{Department of Applied Informatics, Kaunas University of Technology, Studentu 50-415, Kaunas LT-51368, Lithuania}
\date{\today}

\begin{abstract}

We study two coupled systems, one playing the role of the driver system and the other one of the driven system. The driver system is a time-delayed oscillator, and the driven or response system has a negligible delay.  Since the driver system plays the role of the only external forcing of the driven system, we investigate its influence on the response system amplitude, frequency and the conditions for which it triggers a resonance in the response system output. It results that in some ranges of the coupling value, the stronger the value does not mean the stronger the synchronization, due to  the arise of a resonance. Moreover, coupling means an interchange of information between the driver and the driven system. Thus, a built-in delay should be taken into account. Therefore, we study whether a delayed-nonlinear oscillator can pass along its delay to the entire coupled system and, as a consequence, to model the lag in the interchange of information between the two coupled systems. 
 
\end{abstract}

\maketitle

\section{Introduction}

Over the last years an important research activity has been devoted to the dynamics of coupled and driven nonlinear oscillators.  These systems exhibit complex and rich behaviors due to the interplay of nonlinearity, coupling, and external forcing. The results are of interest for several scientific disciplines such as biomedical sciences \cite{jiruska,mormann,Rulkov}, where coupled oscillators are prevalent in biological systems, such as neurons. Studying their dynamics helps us understanding phenomena like neuronal synchronization, leading to advancements in medical treatments and diagnostics. In some cases, coupled oscillators can synchronize their motion, where their frequencies and phases align \cite{jensen,hramov}. To emphasize the role of the previous ideas, we can mention that coupled chaotic oscillators can be employed to enhance communication security \cite{Koronovskii,Naderi}. On the other hand, they are crucial in understanding and controlling vibrations in mechanical systems \cite{defoort}. Furthermore, they have applications in networked systems, ranging from improving the performance of communication networks to understanding the behavior of interconnected systems \cite{delellis,zhang}.Among the fields of interest, we can also mention electronics \cite{Yao} and mechanical engineering \cite{sujith}.

Coupled and driven systems can be modeled as dynamical systems in which one of their parameters is the dynamical variable that comes from another dynamical system, through a coupling mechanism \cite{Pecora1,Pecora2,Pecora3}. Usually, we refer as {\it the driver system} the source of the driving signal, and as {\it the response system} the driven signal. Consequently, the driver system sends a signal to the response system altering its behavior according to the received input. The synchronization of the dynamics of the response system with respect to the driver system~\cite{Boccaletti}, constitutes a relevant effect observed when an oscillator is driven by another oscillator.  Two main cases can be distinguished. When the two oscillators are identical or nearly identical, identical synchronization \cite{Pecora4} may be observed. However, when they are different, generalized synchronization \cite{Rulkov:1995, kocarev:1996} is expected. The driver signal can be either periodic or aperiodic \cite{Ding}. Among the various coupling mechanisms, such as  the replacement method or the subsystem decomposition \cite{Pecora1,Pecora2,Pecora3}, we have chosen to use {\it the continuous control} that was discussed in \cite{Ding,Kapitaniak}. The implementation of the coupling mechanism is carried out by introducing in the response system a square matrix, whose elements are constant, multiplied by the vector of the difference between the dynamical variable of the driver system and the dynamical variable of the response system. 

Nevertheless, the study of the synchronization is not the main goal of this article. In fact, we study the case in which a delayed nonlinear oscillator is the only driver, through the coupling mechanism, of another nonlinear oscillator without delay. Therefore, we analyze the driver system as the only external forcing acting on the response system. We have chosen the continuous control \cite{Ding,Kapitaniak} as a coupling mechanism because we think that it models better the implementation of an external forcing into the system. As a matter of fact, the coupling matrix becomes constant  playing the role of the forcing amplitude, and the time delay determines the frequency of the forcing. In the other above-mentioned methods one or more variables of the driver system are substituted directly into the response system, without the possibility to change the strength of the coupling. This means that there is nothing playing the role of the forcing amplitude. Thus, with that implementation we investigate the typical effects of an external forcing, here affected by delay, on the oscillations amplitude and frequency of a given dynamical system. Moreover, the forcing generates at the right frequency and for the right amplitude the appearance of a resonance. Some applications of implementing delay in an external forcing or control are discussed in \cite{Gu,Sayed1,Sayed2}

Although the synchronization of the two systems is not the main goal of this article, a subsidiary objective can be pondered as a consequence of it. In fact, through the coupling mechanism the driver system transfers to the response system some of its features and here we focus on the delay transmission. We want to emphasize that the synchronization achieved here cannot be identical because the space dimensions of the two systems are different, being one infinite dimensional for the delayed oscillator and the other one finite dimensional. There is a reason to study the conditions of the delay transmission. The coupling, as we wrote before, is an interchange of information from the driver system to the response system and the speed of this interchange is finite, so a built-in delay needs to be considered. Therefore, due to this intrinsic delay, the response system is affected showing delay-induced behaviors. Hence, we have decided to study the optimal parameter values that model this delay in the driving signal and transfer it from the driver system oscillator to the entire coupled system. The result is that, for coupling values that do not trigger the resonance, the response system acts with a similar delay-induced behavior as the driver system, without being a perfect copy. Also, the coupling constant can be used as a control parameter to determine how much delay-induced behavior we want the response system output to show.

The organization of the paper is as follows. In Sec.~\ref{s:model}, we define the model that we have used. We identify the role of the coupling constant in the dynamics of the response system in Sec.~\ref{s:CC}. We discuss in Sec.~\ref{s:CC_tau} the coupling constant and the driver system delay influence on the dynamics of the response system is done. In Sec.~\ref{s:tau}, we analyze some particular values of the coupling constant in function of the driver system delay. Finally, some concluding remarks appear in Sec.~\ref{s:conclucion}.

\section{The model and the continuous control}\label{s:model}

The unidirectional coupling can be summarized in a simple way. We can define an autonomous nonlinear dynamical systems as the driver system and its dynamical state given by a vector $\mathbf{x_1}\in\mathbb{R}^n$ of $n$ scalar variables.
The system dynamics is governed by a set of $n$  nonlinear differential equations $\mathbf{\dot{x}_1}=\mathbf{F}(\mathbf{x_1})$. Then, another nonlinear dynamical system is considered as the response system, whose dynamical
state is given by a vector $\mathbf{x_2}\in\mathbb{R}^d$. The differential equations of this second system are $\mathbf{\dot{x}_2}=\mathbf{G}(\mathbf{x_2})$. When the unidirectional drive is established, the response system becomes $\mathbf{\dot{x}_2}=\mathbf{G}(\mathbf{x_1},\mathbf{x_2})$. 
The continuous control scheme provides a simple form of unidirectional coupling:
\begin{equation}
    \mathbf{G}(\mathbf{x_1},\mathbf{x_2})=\mathbf{G}(\mathbf{x_2})+\mathbf{C}\cdot(\mathbf{x_1}-\mathbf{x_2}),
\end{equation}
where $\mathbf{C}$ is a square matrix of dimension $n$ whose elements are constants. This matrix is multiplied by the vector of differences $(\mathbf{x_1}-\mathbf{x_2})$. The numerical values of the constants inside $\mathbf{C}$ measure the strength of the coupling for each forcing signal, which may be constructed from one, or more, of all the variables of the drive.

\begin{figure}[htbp]
  \centering
   \includegraphics[width=10.0cm,clip=true]{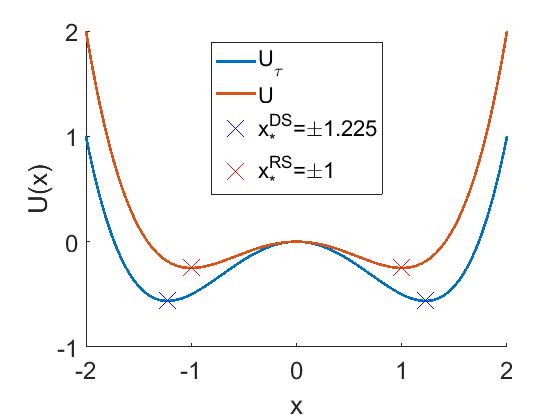}
   \caption{The double-well potentials and stable fixed points of the Duffing oscillator (in red) and of the delayed-Duffing oscillator (in blue).}
\label{fig:1}
\end{figure}

Here, we have decided to study the output of a Duffing oscillator, as the response system, when it is driven by a time-delayed Duffing oscillator,  as the driver system, following:
\begin{align}\label{eq:1}
  Driver&\rightarrow \quad \frac{d^2x_1}{dt^2}+\mu\frac{dx_1}{dt}+\gamma x_1(t-\tau)+\alpha x_1(1-x_1^2)=0\\
  Response&\rightarrow \quad \frac{d^2x_2}{dt^2}+\mu\frac{dx_2}{dt}+\alpha x_2(1-x_2^2)=C(x_1-x_2), 
\end{align}
where we have fixed the parameters $\mu=0.01$, $\alpha=-1$ and $\gamma=-0.5$. The parameter $C$ is the coupling constant, which is the only nonzero element of the continuous control scheme coupling matrix~\cite{Ding,Kapitaniak} that measures the strength of the coupling for the forcing signal and plays the role of the the external forcing amplitude, i.e., the time-delayed Duffing oscillator. The  dissipation $\mu$ is kept small in order to better appreciate the effects of the variation of the driver system delay $\tau$ and of the coupling constant $C$ on the dynamics of the response system. The history functions of the driver system are $u_0=v_0=1$, and the initial conditions of the response system are $x_0=y_0=0.5$.  We expect that our conclusions are of general validity and not specific for the considered boundary conditions.
The potentials 

\begin{align}\label{eq:Pot}
  Driver&\rightarrow \quad \frac{\alpha x^2}{2}+\frac{\alpha x^4}{4}-\frac{\gamma x^2}{2}\\
  Response&\rightarrow \quad \frac{\alpha x^2}{2}+\frac{\alpha x^4}{4} 
\end{align}

and the fixed points of both systems are shown in Fig.~\ref{fig:1}. The unstable fixed point $x_0=0$ is the same for the two potentials, while the stable fixed points of the driver system are
\begin{equation}\label{eq:2}
x^{DS}_{*}=\pm\sqrt{\frac{\alpha+\gamma}{\alpha}}=\pm1.225,
\end{equation}
and the stable fixed points of the response system are
\begin{equation}\label{eq:3}
x^{RS}_{*}=\pm\sqrt{\frac{\alpha}{\alpha}}=\pm1.
\end{equation}
Moreover, following~\cite{Coccolo}, we perform the linear stability analysis for the fixed points. The characteristic equation of the linearized system is
\begin{equation}
   \label{eq:char} \lambda^2+\mu\lambda+\alpha(1+3(x^{DS}_*)^2)+\gamma e^{\lambda\tau}=0.
\end{equation}
Then, we take $\lambda=\rho+i\omega$ as the eigenvalue associated with the equilibrium points. The critical stability curve can be found by fixing $\rho=0$. Hence, we substitute $\lambda=i\omega$ in the last equation and separate the real and imaginary parts obtaining the equations
\begin{align}
    &\omega^2-\alpha(1+3(x^{DS}_*)^2)=\gamma\cos{\omega\tau}\\
    &\mu\omega=\gamma\sin{\omega\tau}.
\end{align}
    After squaring and adding both equations we obtain
    \begin{equation}
   \label{eq:sqare} (\omega^2-\alpha(1+3(x^{DS}_*)^2))^2+(\mu\omega)^2=\gamma^2.
\end{equation}

\begin{figure}[htbp]
  \centering
   \includegraphics[width=15.0cm,clip=true]{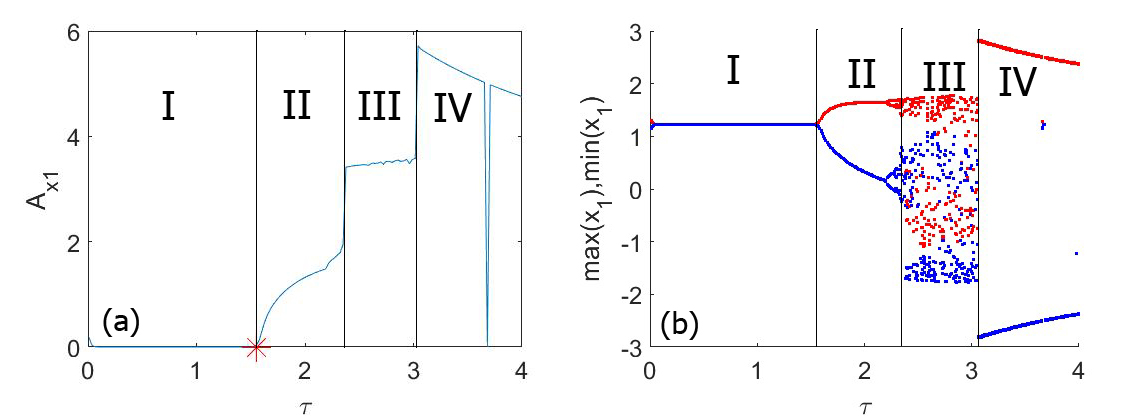}
   \caption{The figures show the oscillations amplitude $A_{x_1}$ (a) and the maxima and minima diagram (b) of the driver system. We can appreciate the oscillations amplitudes (a) and the oscillator behaviors (b) in all the $\tau$ regions of the driver system. The history functions for the time-delayed Duffing oscillator are constant $(u_{0},v_{0})=(1,1)$. The red asterisk in panel (a) is  the value of $\tau$ predicted,  through the Eqs.(\ref{eq:char} - \ref{eq:tau}), by the stability analysis at which the fixed points undergo a change of stability. }
\label{fig:2}
\end{figure}

 Then, substituting the parameter values $\alpha=-1,\gamma=-0.5,\mu=0.01$, we can find four solutions, among which one is $\omega=2.0004$ giving $\tau=1.5505$ as the solution of the equation
    \begin{equation}
     \label{eq:tau}   \tau=\frac{\arccos{((\omega^2-\alpha(1+3(x^{DS}_*)^2)/\gamma)}}{\omega}.
    \end{equation}
 The $\tau$ value just computed is where the fixed points lose stability and is shown as the red asterisk in Fig.~\ref{fig:2}(a).
    
\begin{figure}[htbp]
  \centering
   \includegraphics[width=14.0cm,clip=true]{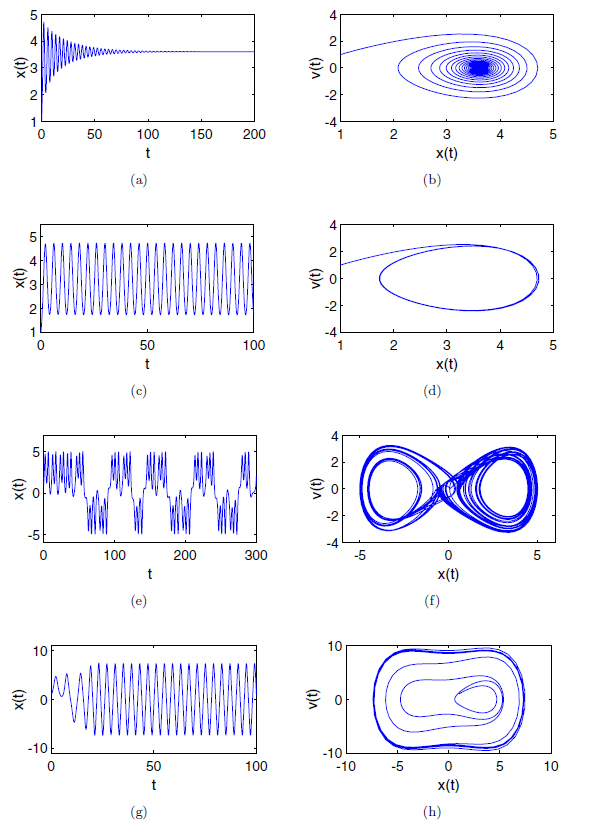}
   \caption{The figure shows the oscillations of the driver system in the stable regions defined in Fig.~\ref{fig:2}. Panels (a) and (b) shows the driver system $x$ oscillations and the orbit in the phase space for $\tau\in$ region I, respectively. Panels (c) and (d) for $\tau\in$ region II. Panels (e) and (f) for $\tau\in$ region III. Panels (g) and (h) for $\tau\in$ region IV. }
\label{fig:2b}
\end{figure}

Then, as already reported in \cite{Coccolo, Cantisan}, the unforced time-delayed Duffing oscillator undergoes various bifurcations while $\tau$ changes. We show such behaviors in Figs.~\ref{fig:2}(a) and (b), where we plot the oscillations amplitude and a diagram showing maxima and minima oscillations amplitudes, respectively. This diagram has been plotted by representing on the figure the maxima and minima of the last $5$ periods of the oscillations for each $\tau$ value. Four regions are discernible in the figures. The first one (\textbf{I}) for $\tau<1.53$ in which the oscillations converge to the fixed point. The second one (\textbf{II}),  $1.53<\tau<2.35$ where the oscillations are sustained and confined to one of the wells. The third one (\textbf{III}),  $2.35<\tau<3.05$ the amplitude has jumped to a value bigger than the width of the well, what it means that the trajectories move from one well to another one, and also they are aperiodic. The last one (\textbf{IV}), for $\tau>3.05$ where the trajectories are no longer confined to either of the wells and oscillations are periodic. The result is a limit cycle in phase space that spans both wells. All the simulations of the manuscript have been carried out with the DDE tools of Matlab and checked with a fourth-order Runge-Kutta integrator for the non-delayed case, with an integration step of $0.01$. These behaviors of the driver system are depicted in Fig.~\ref{fig:2b}.

\section{The role of the coupling constant}\label{s:CC}

\begin{figure}[htbp]
  \centering
   \includegraphics[width=15.0cm,clip=true]{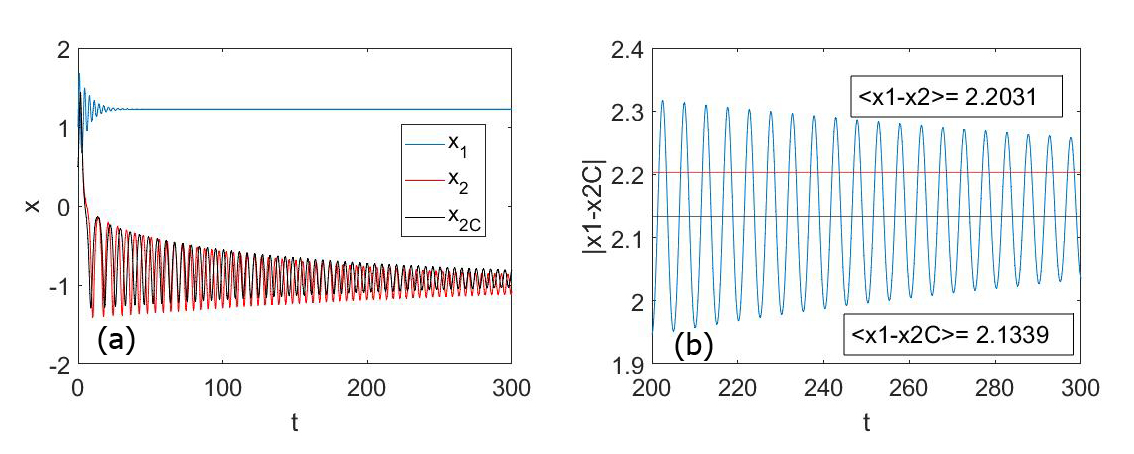}
   \caption{Panel (a) shows the driver system and the response system $x$-time series at fixed $\tau=1$. In blue the driver, in red the Duffing with $C=0$ and in black the response system with $C=0.06$. Panel (b) shows the distance between the driver system and the response system $x$-time series $|x_1(t)-x_{2C}(t)|$ (the blue oscillations) and the mean distances $<x_1-x_2>$ and $<x_{1}-x_{2C}>$, the red and black lines, respectively.  }
\label{fig:3}
\end{figure}

In Fig.~\ref{fig:3}(a), we show the $x$-time series of the two systems when $C=0$, i.e., without coupling, and with $\tau=1$ and $\mu=0.01$. In the figure, we can see in blue the driver system $x$-time series oscillations, that from now on we call $x_1$, while in red the response system without coupling, that we call $x_2$. We can appreciate how the two oscillators tend to their specific fixed points independently one from the other. 
Once we have seen how the two oscillators behave independently, from this point on we switch on the coupling constant so that the time-delayed Duffing starts to drive the Duffing oscillator without delay. The effects are shown in Fig.~\ref{fig:3}, in which in black we represent the response system $x-$time series for $C=0.06$, from now on $x_{2C}$.  The coupling value has been chosen for explanatory purposes. In fact, it is the value for which it is possible to appreciate that the oscillations of $x_{2C}$ are slightly displaced up towards the fixed point of the driver system, although they have not yet jumped into the other well.  In Fig.~\ref{fig:3}(b), we show (the curve in blue) the absolute value of the asymptotic distance between $x_1$ and $x_{2C}$, $|x_1-x_{2C}|$, for $t>200$. Also, it is shown that the mean distance between  the $x$-time series of the coupled response system and the $x$-time series of the driver system, black line, $<x_1-x_{2C}>=2.1339$, is smaller than the mean distance of the $x$-time series of the uncoupled response system and the $x$-time series of the driver system, red line,$<x_1-x_{2}>=2.2031$. This measures the level of synchronization between the two oscillators, following the standard definition of synchronization \cite{Miranda},
\begin{equation}
\label{eq:5} \lim_{t\to\infty}|x_1(t)-x_{2C}(t)|\rightarrow0, 
\end{equation}
stating that if the mean of the asymptotic distance between the solutions of the two oscillators goes to zero, the two oscillators are synchronized. From now on, when we state that one case is more synchronized than another, it means that this definition has been used.
To obtain the mean of the asymptotic distance, we have computed the absolute value of the difference between the last third part of the $x$-series of the two systems and then its mean value.

\begin{figure}[htbp]
  \centering
   \includegraphics[width=15.0cm,clip=true]{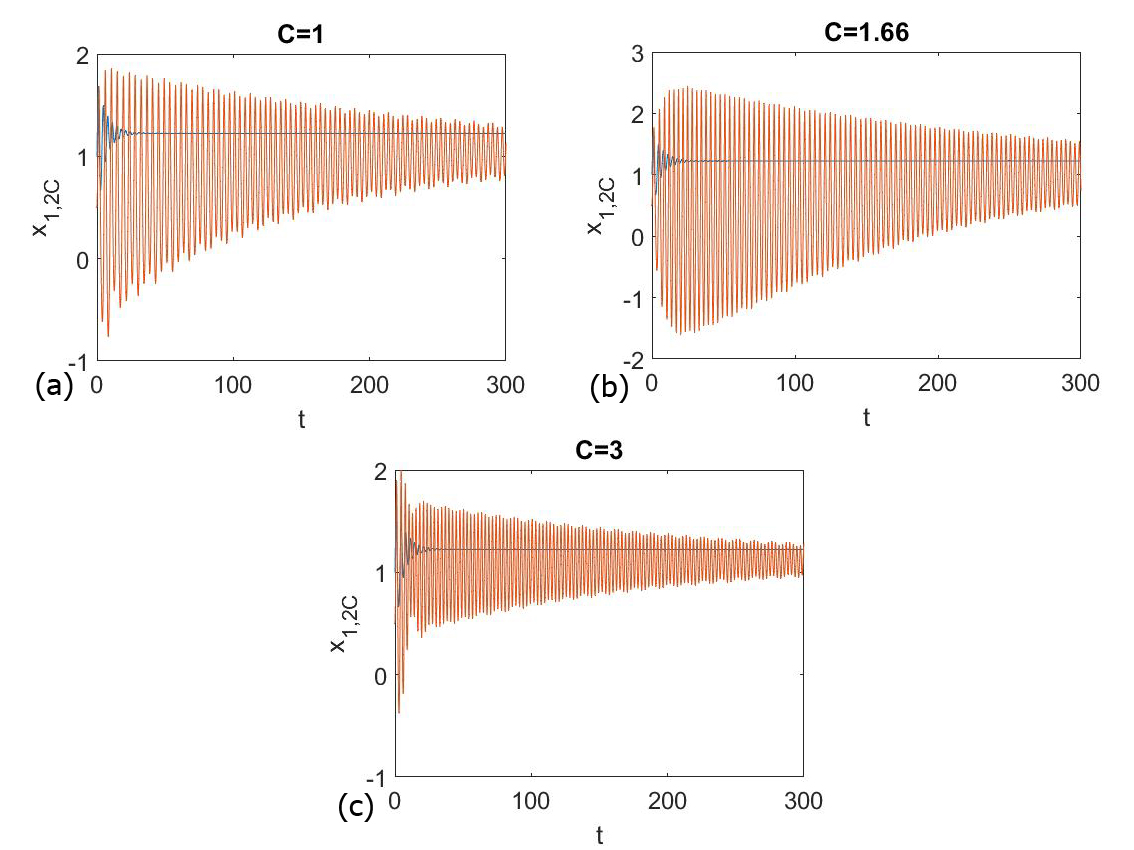}
   \caption{The figure shows the effect for growing coupling constant on the dynamics of the response system, $C=1, C=1.66$ and $C=3$ at fixed $\tau=1$. The panels show the driver system $x$-time series, $x_1$ in blue, and the response system  $x_{2C}$ in red. In the panels it is appreciable that the oscillations amplitude of the response system is larger in the $C=1.66$ case. Moreover, counter intuitively the synchronization of the two system is better at $C=1$ than at $C=1.66$. }
\label{fig:4}
\end{figure}

Then, we set larger values of $C$ for $\tau=1$, as shown in Fig.~\ref{fig:4}. We can see that for $C=1$, Fig.~\ref{fig:4}(a), the response system $x$-time series jumps into the other well, its main distance from the driver system becomes significantly smaller $<x_1-x_{2C}>=0.2528$ than in Fig.~\ref{fig:3}(b). Certainly, the main distance between $x_1$ and $x_2$ does not change. Contrary to the expectations when the $C$ value increases at $C=1.66$, Fig.~\ref{fig:4}(b), the response system $x$-time series asymptotic oscillations grow larger and the mean distance $<x_1-x_{2C}>=0.4733$ increases. Then, if the coupling constant increases further at $C=3$,  Fig.~\ref{fig:4}(c), the asymptotic oscillations amplitude decreases and also the mean distance $<x_1-x_{2C}>=0.1460$, as it is supposed to be.   In other words, comparing the Figs.~\ref{fig:4}(a)-(c), the oscillations amplitude of the response system grows, reaches a maximum and decreases in function of the coupling constant. This effect of the coupling constant is comparable to the effect of an external forcing amplitude that induces a resonance.

\begin{figure}[htbp]
  \centering
   \includegraphics[width=14cm,clip=true]{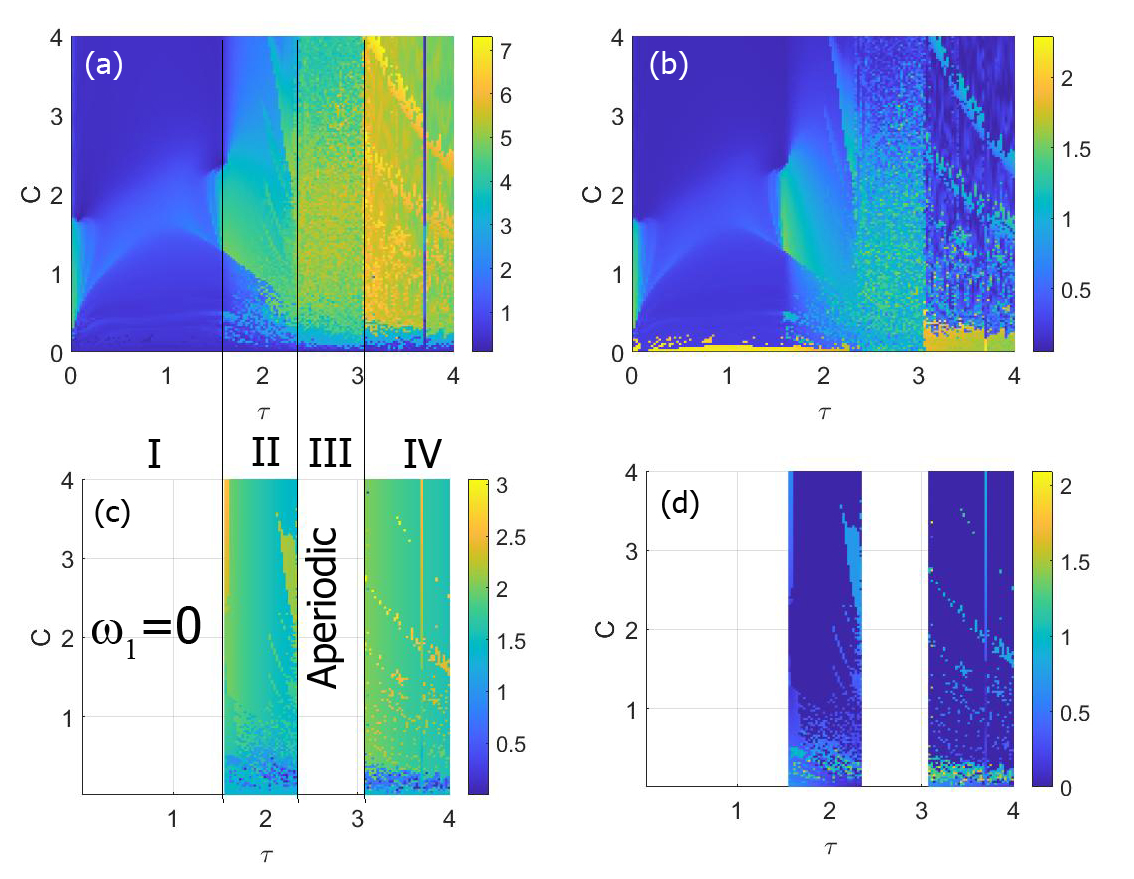}
   \caption{The panels show (a) the gradient of the amplitude of the response system oscillator $A_{x_{2C}}$.  Then, (b) the gradient of the mean distance between the driver system and the response system. Finally (c), the gradient of the response system frequency $\omega_{2}$, and (d) the gradient of the difference between the frequency of the driver system and the response system, $|\omega_1-\omega_{2}|$. The frequencies have been calculated using the fast Fourier transform. All figures show gradient plots in function of the coupling constant $C$ and the time delay $\tau$.   In panels (c) and (d), only the regions II and IV are represented, because the driver system frequency in region I is zero and region III is aperiodic, making the comparison with the response system meaningless.}
\label{fig:5}
\end{figure}

\section{The combined effect of the coupling constant and the delay}\label{s:CC_tau}

Now to start the analysis of the mentioned coupling constant effect and its interaction with the delay $\tau$, we plot Fig.~\ref{fig:5}. Here, we can find the following gradient plots. In Fig.~\ref{fig:5}(a), the oscillations amplitude of the response system, $A_{x_{2C}}$. In Fig.~\ref{fig:5}(b) the mean distance between the driver system and the response system asymptotic behaviors, $<|x_1(t)-x_{2C}(t)|>, t>200$. In Fig.~\ref{fig:5}(c) the response system oscillation frequency, $\omega_{2}$, for the driver system regions with periodic oscillations. All throughout the manuscript the frequencies of the response system and of the driver are all calculated using the fast Fourier transform.  In Fig.~\ref{fig:5}(d) the difference between the oscillations frequencies of the driver system and response system, $|\omega_1-\omega_{2}|$.  All the gradient plots are in function of the delay $\tau$ and the coupling constant $C$. We can see in Fig.~\ref{fig:5}(a) that the oscillations amplitude of the response system grows with the delay $\tau$ of the driver system, similarly to the driver system itself, as shown in Fig.~\ref{fig:2}.  On the other hand, the oscillations amplitude grows and changes for some range of the coupling constant $C$ that varies in every region. Then, Fig.~\ref{fig:5}(b) shows the level of synchronization between the two oscillators. In fact, every point in the figure represents the mean distance of the asymptotic behaviors of the oscillators $x$-series. So that, the smaller the mean distance between the two oscillators $x$-series, the higher the level of synchronization. Counter intuitively, it is not clear all along the panel that the larger the coupling constant the better the synchronization. In fact, the synchronization is better in some regions like region I, for $C = 1$ and $\tau = 1$, than for the case $C = 1.66$ for the same $\tau$ value. To better understand the previous figure, we analyze the panels region by region and we study a particular case for a fixed $\tau$ and varying $C$ in the interesting cases.

\subsection{The coupling constant effect}
\begin{description}[style=unboxed,leftmargin=0cm]

\item[In region I] the oscillation amplitudes, Fig.~\ref{fig:5}(a), are smaller with respect to the other regions. The exception being the zone of lower $\tau$ values and $0\lesssim C\lesssim2$. Then, an area of relatively higher amplitude is visible in the middle of the region.  In fact, we can appreciate a higher amplitude “tubular structure” that cross the region along the $\tau$ values, but only for certain $C$ values. This peak resembles a resonance peak. As a matter of fact, the response system oscillations are small outside “the tube” but are much larger inside it, as a consequence of the the driver system. As expected for a resonance, there is a dependence on the amplitude of the external forcing, in this case $C$. However, there is also a dependence on the frequency of the external forcing, here the delay $\tau$.  In Fig.~\ref{fig:5}(b), we can notice a high difference between the two oscillators $x-$series for $C$ values outside  and  inside the resonance area. So that, smaller $C$ values can synchronize the two oscillators more than larger $C$ values, located inside the resonance area, as in the already discussed example of $\tau=1$ in Fig.~\ref{fig:4}. 

Now, we study a particular case, for which we consider for a given value of $\tau=1$ a range of $C$ values shown in Fig.~\ref{fig:5}. This $\tau$ value lies in the region I depicted in Fig.~\ref{fig:2}. Thus, we plot in Fig.~\ref{fig:6}, the asymptotic oscillations amplitude and the asymptotic behavior of the response system oscillations as a maxima-minima diagram. All of this for the coupling constant values $0.01<C<4$. In particular, Fig.~\ref{fig:6}(a) shows the amplitude $A_{x_{2C}}$ in which we can recognize a bell shaped curve that reminds of a resonance, with its maximum at $C=1.66$ while the red line is the amplitude $A_{x_{2}}\simeq 0.57$ of the response system with $C=0$. Finally, the black line is the amplitude of the driver system, $A_{x_1}\simeq 0$, which oscillations fall into the fixed point. This peak is a section of the $\tau$ space spanning peak seen in Fig.~\ref{fig:5}. The peak appearance is recognizable in the maxima-minima diagram, Fig.~\ref{fig:6}(b). These last figures have been portrayed by plotting the maxima and minima of the oscillations last $5$ periods, to show the oscillators asymptotic behaviors. In the figure, also appear as a straight black line the driver system asymptotic behavior that is constant because the change of $C$ does not affect it. Also, in Fig.~\ref{fig:6}(a) we can spot a little peak around $C=0.444$ that matches with a change in the maxima-minima diagram of the response system in Fig.~\ref{fig:6}(b). This little peak is related  with the small filiform zone around $C\thickapprox0.444$ that crosses all the region I along the $\tau$ axis in Fig.~\ref{fig:5}(a). All that has been written shows how the coupling constant is introducing a perturbation into the response system that, in addition to driving it, also forces the system. This triggers a resonance between the driver system and the oscillations of the response system. 

\item[In region II] we can see that a zone of higher amplitudes, Fig.~\ref{fig:5}(a), starts for values of $C\thickapprox1$ and is connected to the high amplitude area of the region I. In fact, we can recognize that the resonance area of region I continues in the region II. These amplitude maxima spread for more $C$ values until, in the right of the figure, it reaches a wide zone before merging with the aperiodic region III. Then, we can distinguish a variety of peaks in the response system amplitude, most of them for $C\lesssim 1$. These peaks seem generated by an erratic behavior in function of the coupling values of the response system in response of the driver system. We call them {\it adjustment peaks}, because the response system is adjusting its behavior to the driver system before reaching the higher amplitude trend of the higher $C$ values. In these peaks zone, a little variation in the coupling constant can be determinant for the response system to change its asymptotic behavior into the driver system well. The name adjustment peaks comes from the fact that the coupling constant $C$ for the peaks cannot overdue the dynamics of the response system, but it is able to stretch the orbits at a larger amplitude. On the other hand, a slightly bigger or smaller value of $C$  can drive the response system into the driver system well, so that the response system amplitude oscillates in function of $C$. 
In region I and II, the gradient plot of the oscillations difference $<|x_1(t)-x_{2C}(t)|>$, Fig.~\ref{fig:5}(b), shows that to obtain a better synchronization of the two oscillator values of $C$ outside the resonance area should be chosen. In Figs.~\ref{fig:5}(c) and (d), we can appreciate that the response system oscillation frequency $\omega_2$ grows and the difference between the two oscillators frequencies, $|\omega_1-\omega_2|$ decreases while $C$ grows up. However, Although, the difference grows at $\tau\thickapprox1.5$ and $C\thickapprox4$, because we are close to the aperiodic region III and some fluctuation can start. This behavior is not visible by the difference in amplitudes, Fig.~\ref{fig:5}(b), so the effect is just on the frequency.
Here, we repeat what we have done before, we study a particular case, here $\tau=2$. In fact, in Figs.~\ref{fig:7}(a) and (b) we can see how the amplitude of the response system vary in function of the coupling constant $C$. The results can be visualized in the maxima-minima diagram, as shown in Fig.~\ref{fig:7}(b). As before, the black straight lines that are the driver asymptotic behaviors have been plotted for comparison with the response systems in Fig.~\ref{fig:7}(b) and the behavior of the driver system does not change for different values of $C$. In this set of figures, we can distinguish a variety of peaks in the response system amplitude, most of them for $C\leq 1$. These peaks are the adjustment peaks that we described before. The peak at $C=1$ is different because it starts the general tendency that goes further for bigger values of $C$, i.e., the amplitude decreases in order to adjust to the amplitude of the driver system. However, for $C>1$ some little peaks call our attention since they disrupt the general tendency of the curve, in particular the one for $C=3$. This can be generated by a resonance between the external forcing (the driver system) and the system (the response system). Finally, in Fig.~\ref{fig:7}(c), we show the change in the frequency $\omega_2$ of the response system, in function of $C$. Here, we can appreciate that, initially, $\omega_2=1.3684$, then its value oscillates until it becomes $\omega_2=\omega_1=1.5881$, for $C>1$. We have to take into account that the frequency plot has been only analyzed in this case and for $\tau$ values in the region IV, because just in these two regions the oscillations of the driver system are periodic.

\item[In region III]  the $\tau$ values are in the aperiodic region, defined in Fig.~\ref{fig:2}, so the response system behaviors become aperiodic when the coupling value exceeds $C\thickapprox0.06$, see Fig.~\ref{fig:5}(a) and Fig.~\ref{fig:5} (b).

\item[In region IV] we can find high oscillation amplitudes for both oscillators, see Fig.~\ref{fig:5}(a).  Here, we can see that high oscillations amplitude grows for $C\gtrsim0.06$. Besides, the driver system shows a minimum in the oscillations amplitude at $\tau=3.68$, and the same behavior appears in the response system for almost all the $C$ values. In Fig.~\ref{fig:5}(b) we can see how the difference between the oscillators $x-$series is maximum for very small values of $C$, but then it starts to oscillate, growing and shrinking when peaks in amplitude show up and finally becomes small for $C\thickapprox4$. As a result, the response system behaviors adjust to the driver system behavior in a complicated way inside the chosen range of $C$ values. In fact, there are a lot of fluctuations in the response system amplitude,  Fig.~\ref{fig:5}(a), and in the mean difference $<|x_1(t)-x_{2C}(t)|>$, Fig.~\ref{fig:5}(b), except for $\tau\thickapprox4$ and $C\thickapprox 4$. Similar oscillating behaviors can be found in the frequency plots, Figs.~\ref{fig:5}(c) and (d). In the last phrase we do not take into account values of $\tau$ near the region III because the aperiodic oscillations influence is still noticeable.

\end{description}

\subsection{The $\tau$ effect}

We know that, in a time-delayed system, the variation of the delay $\tau$ is responsible for the modification of the frequency of the delay-induced oscillations \cite{Cantisan}, when they exist. So, an interesting question is: how do the resonance peaks spotted in the previous figures behave when $\tau$ varies?  To answer this question, we study the oscillations of the response system while the external forcing frequency, the driver system delay, changes.
Therefore,  we analyze Fig.~\ref{fig:5} along the $\tau$ axis.

\begin{figure}[htbp]
  \centering
   \includegraphics[width=14.9cm,clip=true]{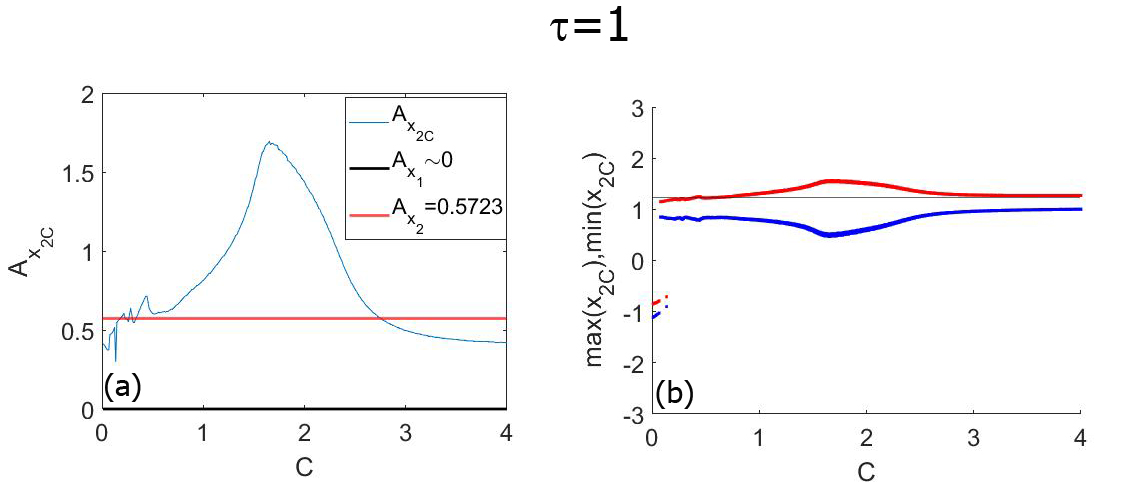}
   \caption{The panels show a slice of the previous Fig.~\ref{fig:5} for varying values of $C$ and fixed $\tau=1$. In particular (a) it portraits the amplitude of the oscillators, being $A_{x_{2C}}$ the amplitude of the response system with the nonzero coupling term, $A_{x_{2}}$ the amplitude of the response system without the coupling term and $A_{x_1}$ the amplitude of the driver system. Then, (b) the maxima and minima diagram of the asymptotic behavior for the driver (the black straight line) and the response system. It is interesting to observe the amplitude peaks in panel (a): the smaller one for $C=0.444$ and the larger one for $C=1.66$  that suggests a resonance induced by the coupling of the two oscillators. In panel (b), the maxima and minima of the driver system overlap since the oscillator goes to the fixed point. It is also shown that the oscillations of the response system tend to the driver system fixed point (straight black line), for bigger values of $C=0.06$ and values outside the resonance area.  }
\label{fig:6}
\end{figure}

\begin{figure}[htbp]
  \centering
   \includegraphics[width=13.5cm,clip=true]{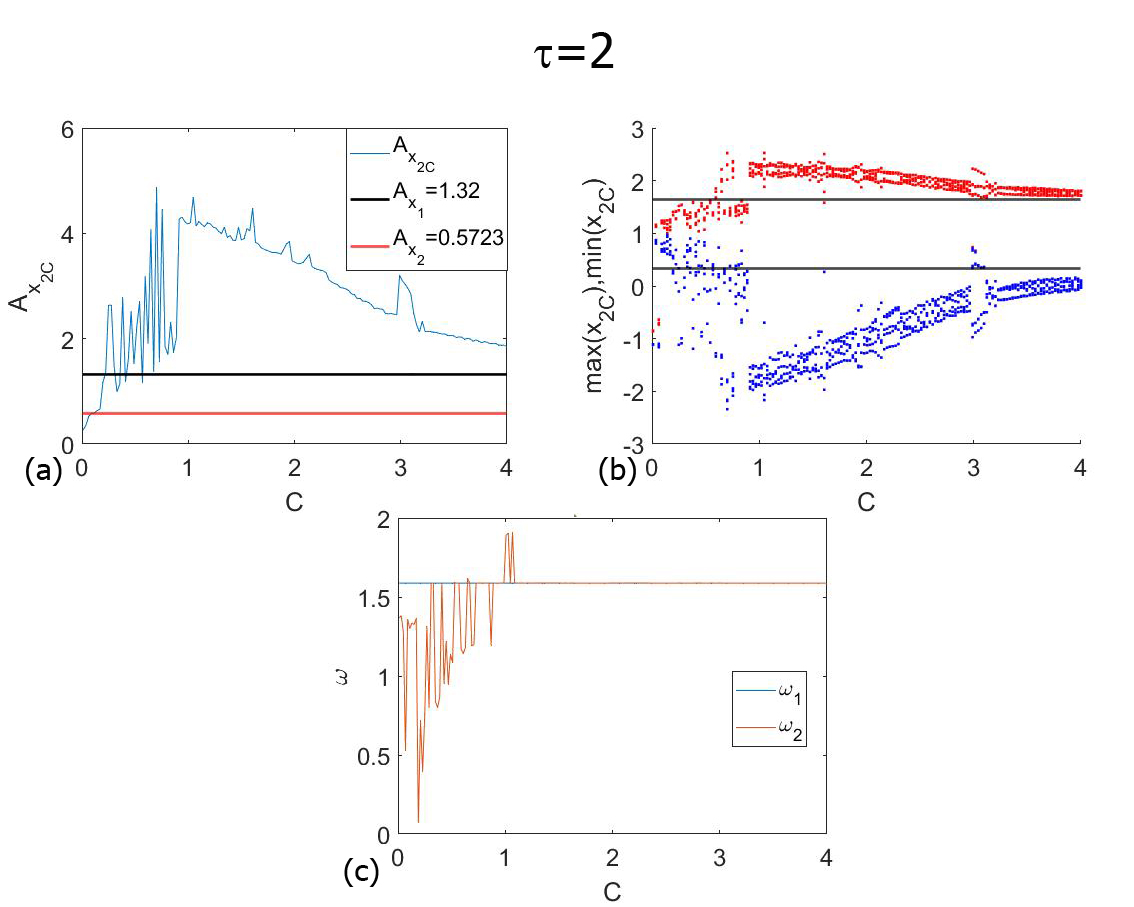}
   \caption{The panels show (a) the amplitude of the oscillators, being $A_{x_{2C}}$ the amplitude of the response system with the nonzero coupling term, $A_{x_{2}}$ the amplitude of the response system without the coupling term and $A_{x_1}$ the amplitude of the driver system. Then, (b) the maxima and minima diagram of the asymptotic behavior for the driver (black straight lines)  and the response system for varying values of $C$ and fixed $\tau=2$, respectively. Some interesting peaks appear in the amplitude plot, due to a resonance. Finally, panel (c) shows the $\omega_1$ of the driver system in blue, that does not change when $C$ changes and $\omega_2$ of the response system in red that for $C=0$ is $\omega_2=1.3684$ and when the coupling constant is larger than $C\simeq1.2$ it becomes $\omega_2=\omega_1$.}
\label{fig:7}
\end{figure}

\begin{figure}[htbp]
  \centering
   \includegraphics[width=13.5cm,clip=true]{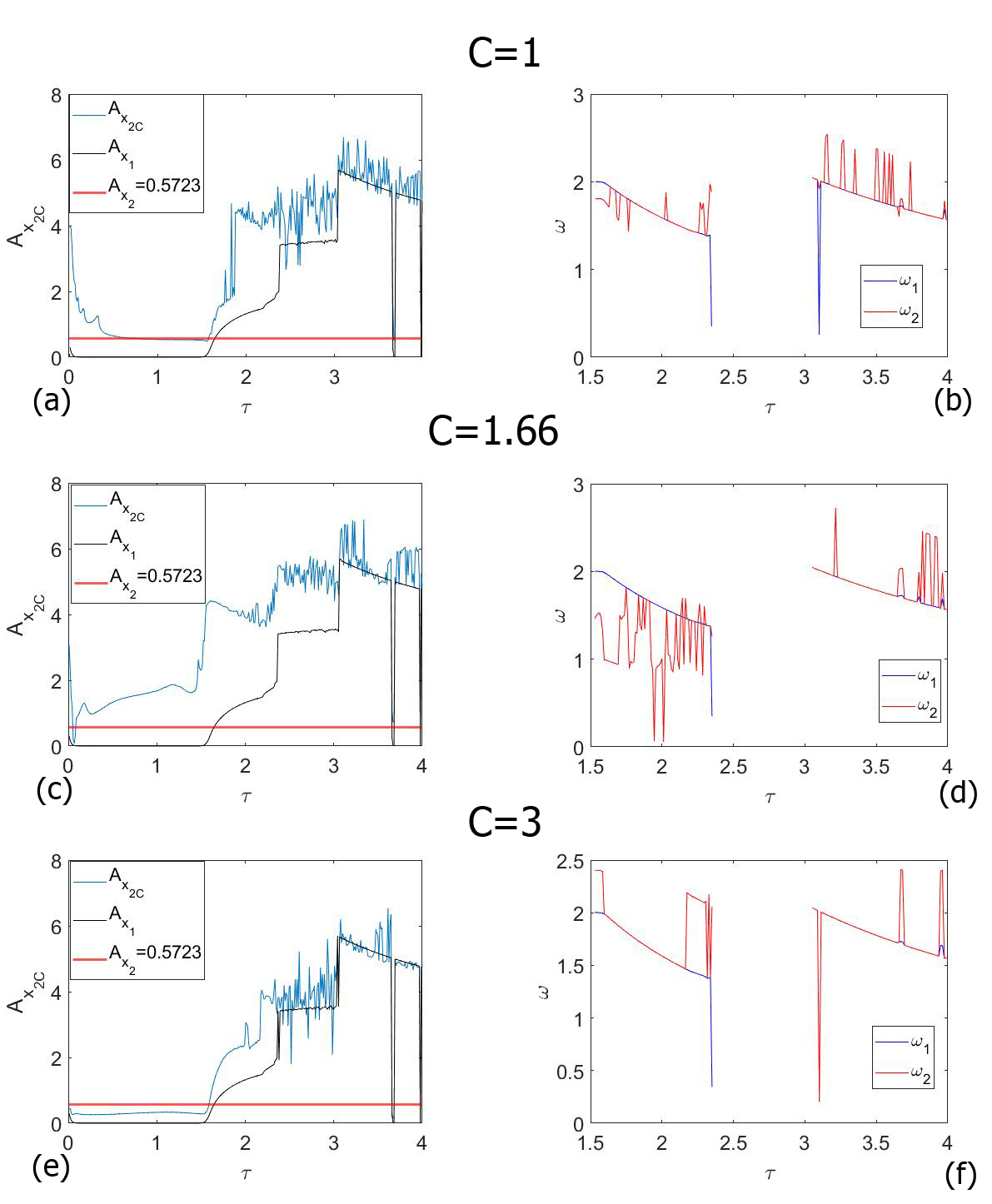}
   \caption{Panels (a), (c) and (e) show the oscillations amplitude and panels (b),(d) and (f) the frequency in the region II and IV for three coupling values. The first one before the peak in Fig.~\ref{fig:6}(a), i.e., $C=1$. The second is the value that gives the top of the peak, $C=1.66$. The third one beyond the peak, $C=3$. Interestingly in panel (a) we can see that in the region I  for $C=1$, the response system follows the driver better than the second case, panel(c), although the second coupling constant value is larger than the first one. The resonance peaks already mentioned are recognizable in Region I, II and III for $C=1.66$ and in region II for $C=1$. Finally, it is interesting that in Region II the frequencies follows the $\omega_1$ better than in all the other cases.}
\label{fig:8}
\end{figure}

We start this second part for a coupling constant $\mathbf{C=1}$ at the beginning of the Fig.~\ref{fig:6}(a) resonance peak and by extension on the border of the tubular structure described in Fig.~\ref{fig:5}(a). This coupling value gives a particular peak in the amplitude plot for $\tau=2$,  Fig.~\ref{fig:7}(a). Now by analyzing  Fig.~\ref{fig:8}(a), we can see that the response system amplitude in region II are high, and the driver system is forcing the response system to oscillate at large amplitudes between its fixed point and the one of the driver system. On the other hand, in the region I the oscillations amplitude follows the driver. In fact, we can see that the $\tau$ value for which there is a change in the stability for the response system corresponds with the value predicted by the stability analysis of the driver. Analyzing the frequencies, Fig.~\ref{fig:8}(b), we can see that there is a correspondence of the values in region II. It is worth to mention a down spike in the $\omega_1$ at $\tau\thickapprox3.1$, not appearing in all the $\omega-\tau$ plots, or being different in case that it appears. It seems that so close to the chaotic region III some of its influence remains. In fact, the value of the frequency is still variable for $\tau$ values slightly larger than $3.05$. Important to mention is that this is the coupling value for which the curves of the two frequencies are more similar in region II, except for some values. The most interesting exception is at $\tau=2$, where the oscillations amplitude reaches a maximum, it is also possible to spot a peak for the frequency $\omega_2$,  Fig.~\ref{fig:8}(b).

Keeping on our analysis, bring us to the following value of $\mathbf{C=1.66}$, see Figs.~\ref{fig:8}(c) and (d). Value for which in the  Fig.~\ref{fig:6}(a) a maximum is reached by the response system amplitudes and that falls inside the tubular structure of Fig.~\ref{fig:5}(a). In Fig.~\ref{fig:8}(c), we can appreciate that the oscillations amplitude, in region I, in region II and in region III, are enhanced with respect to the previous value of $C$. Contrary to the expectations, a larger coupling value with respect to the former one does not give a better synchronization of the two systems due to the presence of the resonance. In fact, in the above-mentioned regions, the response system oscillations amplitude, at this value of the coupling constant, reaches its maximum and the $\omega_1$ and $\omega_2$ curves in region II lose coherence, Fig.~\ref{fig:8}(d). In the regions in which the resonance high amplitude are present, the two oscillators are no longer synchronized in both the oscillations amplitude and the frequency.

The last case is for $\mathbf{C=3}$, Figs.~\ref{fig:8}(e) and (f), related with the small peak in Fig.~\ref{fig:7}(a). Here, the coupling constant is big enough to fall outside the resonance tubular structure described in Fig.~\ref{fig:5}(a) and far beyond the values that produce the peak in Fig.~\ref{fig:6}. In fact, the oscillations amplitude are smaller than the previous case and since the coupling constant is large, in the range that we have used, the oscillations amplitude of the response system follows the amplitude of the driver system better than in all the other cases. Also, the $\tau$ value for the response system changes the stability corresponds with the driver system value, see Fig.~\ref{fig:8}(e).  In this case, the driver system and the response system frequencies match almost as well as the $C=1$ case in the region II or better in region IV. Definitely, they match better than in the $C=1.66$ case. So, again, we recognize how the oscillations amplitude and the difference in oscillations between the driver and the response system grow, reach a maximum and decrease in function of the coupling constant and how they vary in function of $\tau$ has also been studied.

\begin{figure}[htbp]
  \centering
   \includegraphics[width=14cm,clip=true]{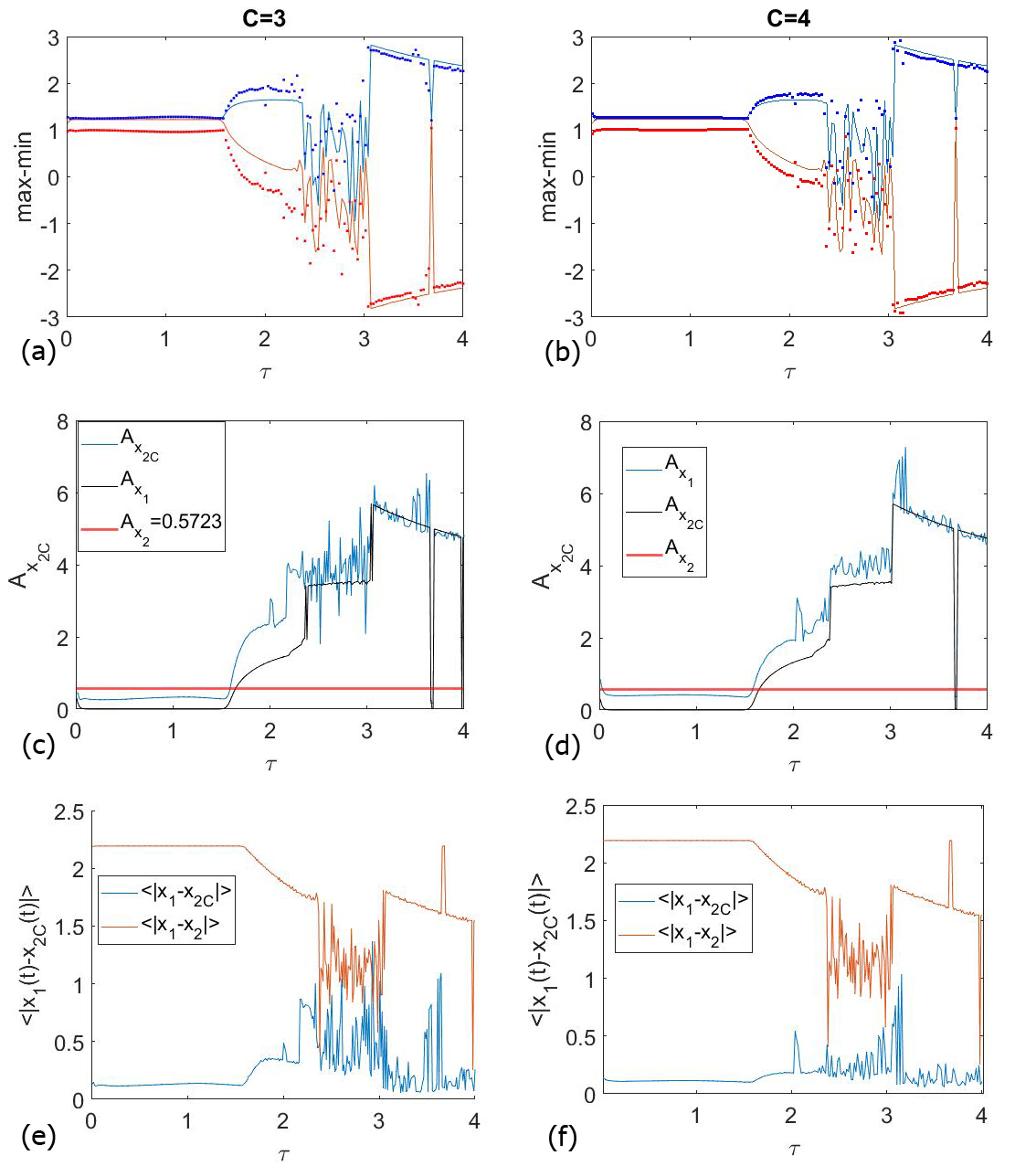}
   \caption{The panels (a) and (b) show, respectively, for $C=3$ and for $C=4$ the mean of the maxima and minima of the last 5 periods of the driver system (solid line) and  of the response system (the dots). Panels (c) and (d) show in black the oscillations amplitude of the driver and in blue of the response system, the red line is the amplitude of the Duffing with $C=0$. Panels (e) and (f) show the differences of the $x-$series of the driver and the Duffing in red, and of the driver and the response system in blue.  For this constant coupling values the response system dynamics adapt to the driver system and the delay-induced oscillations are completely transmitted from the driver system to the response system. Hence, we have modeled the built-in delay of the synchronization just using the delay in the driver system.}
\label{fig:9}
\end{figure}
\FloatBarrier
\section{Delay-induced oscillations in a non-delayed system}\label{s:tau}

In the introduction we wrote about a subsidiary objective of this work. We use the delayed driver system as the only excitation of the driven system and we want to ascertain for which values of the coupling constant its features are better transferred. In particular, coupling means an interchange of information between the driver and the driven system and the speed of this interchange is finite. Thus, an in-built delay should be taken into account. This is the reason to study whether a delayed-nonlinear oscillator can pass along its delay to the entire coupled system. 
In previous sections, we have seen that the best candidates are values of $C$ outside the resonance areas.  This means that through the coupling mechanism the driver system can transfer some of its delay-induced oscillations to the response system in a complicated way that depends on the $\tau$ region. Thus, the response system starts to behave as a delayed oscillator.
To visualize this effect we can focus on Fig.~\ref{fig:8}(e). Then, to analyze deeply the phenomenon we plot  Fig.~\ref{fig:9}. In the panels (a) and (b), we show the mean of the maxima-minima diagram of the last $5$ periods of the orbit in function of $\tau$. The lighter line on the background is the mean of the maxima-minima diagram of the driver system.   In Figs.~\ref{fig:9}(c) and (d) we show the oscillations amplitude of the driver and the response system. In Figs.~\ref{fig:9}(e) and (f) the distance between the driver and the response system $x-$series in blue and the distance between the driver and the Duffing for $C=0$ in red. In these panels we can appreciate a good agreement between the driver and the response system. Looking at those figures we can assure that, for those values of $C$, the delay is completely transmitted, although not perfectly, from the driver system to the response system. In Figs.~\ref{fig:9}(a) and (c), in the region II the effect of the Fig.~\ref{fig:7} peak at $C=3$ is visible as a zone of higher amplitudes. In fact, slightly smaller or larger values of $C$ outside the mentioned peak can guarantee a good enough match with the driver system as the $C=4$ case.
So, we can say trustingly that it is possible for the driver to carry on the synchronization built-in delay just with its own delay. As we saw, the coupling constant value has to fall outside the regions of resonance. If the coupling constant is large enough, i.e., $C>3$ the synchronization of the delay-induced oscillations is not just localized to a particular $\tau$ region or value, as we have seen the case of $C=1$ and $\tau$ region 1, but it is generalized to all the regions. For $C>4$ this effect does not change drastically in comparison with the case $C=4$.
Finally, we also show here that we can use the coupling constant to control which degree of delay-induced oscillations we want to transfer to the response system in order to control the in-built delay due to the coupling.

\section{Conclusions}\label{s:conclucion}

Two coupled systems have been studied. A time-delayed Duffing oscillator as the driver system and a Duffing oscillator without delay as the response system. The driver system plays two roles, the first one as the external forcing of the response system and the second one as the responsible to bring the coupling built-in delay into the response system. 

As regards the first role, the driver system behaving as an external forcing can induce a resonance in the response system as we have seen in the case of $C=1.66$ in the regions I, II and III. Also, in the case of $C=1$ a resonance shows up but just in the region II. Finally, some other resonance peaks pop up, in region I and II, that can be assimilated, in the case of a periodic external forcing, as peaks related with other harmonics. An interesting feature is the adjustment peaks, that give birth to fractal-like zones in the $C-\tau$ gradient plots. In those zones, the amplitude of the response system is highly sensitive to the $C$ value. In fact, for really close values of the coupling constant the response system can fall into the driver system well or have the orbits stretched by the influence of both the driver system and of the response system fixed points.

On the other hand, when the coupling constant takes values outside the resonance area the delay-induced behaviors are better transferred from the driver system to the response system. So that, the response system behaves as a time-delayed system, allowing for the finite velocity in the coupling transmission of information. The best synchronization along all the $\tau$ value is reached at $C=4$, but there are specific cases, as $\tau$ regions or values that reach a good synchronization for smaller $C$ values. As it has been shown, the difference between the asymptotic oscillations of the two systems changes in the parameter set in a complicated way. Finally, the coupling constant can be used as a control parameter that allowing to decide how much of the delay-induced oscillations we want to be transmitted from the driver to the response system, always remembering that identical synchronization is impossible, as explained before. This means that a perfect transmission of the delay from the driver system to the response system is unattainable.

To summarize, both roles unveiled interesting properties in the response system behavior. First of all, the coupling mechanism can trigger a resonance in the driven system oscillations as an external forcing. Second, it is possible to model the delay due to the coupling interchange of information, just with the delay of the driver. Also, we have proved that a previous study can give suggestions to which values of the coupling constant we should use in a specific case. Not always it is possible to apply the stronger coupling value. Also, we can decide, through the coupling constant, how much of the delay-induced oscillation are transferred from one oscillator to the other.

\section{Acknowledgment}
This work has been supported by the Spanish State Research Agency (AEI) and the European Regional Development Fund (ERDF, EU) under Project No.~PID2019-105554GB-I00  (MCIN/AEI/10.13039/501100011033).

%\FloatBarrier


\begin{thebibliography}{99}
%1
\bibitem{jiruska} Jiruska, P., de Curtis, M., Jefferys, J.G.R., Schevon, C.A., Schiff, S.J. \& Schindler, K. Synchronization and desynchronization in epilepsy: controversies and hypotheses.  J. Physiol. 591, 787-797 (2013).
\bibitem{mormann} Mormann, F., Lehnertz, K., David, P., \& Elger, C.E.  Mean phase coherence as a measure for phase synchronization and its application to the EEG of epilepsy patients. Phys. D: Nonlinear Phenom. 144, 358-369 (2000).
\bibitem{Rulkov}Rulkov, N. F.  Regularization of synchronized chaotic bursts. Phys. Rev. Lett. 86, 183 (2001).
\bibitem{jensen} Jensen, R.V. Synchronization of randomly driven nonlinear oscillators. Phys. Rev. E 58, R6907 (1998).
\bibitem{hramov} Hramov, A.E., and Koronovskii, A.A. An approach to chaotic synchronization. Chaos 14, 603-610 (2004).
\bibitem{Koronovskii} Koronovskii, A.A., Moskalenko, O.I., \& Hramov, A.E. On the use of chaotic synchronization for secure communication. Phys.-Usp. 52, 1213 (2009).
\bibitem{Naderi}Naderi, B., and Kheiri, H.  Exponential synchronization of chaotic system and application in secure communication. Optik 127, 2407-2412 (2016).
\bibitem{defoort} Defoort, M., Hentz, S., Shaw, S.W., \& Shoshani, O. Amplitude stabilization in a synchronized nonlinear nanomechanical oscillator. Commun. Phys. 5, 93 (2022).

\bibitem{delellis} DeLellis, P., Di Bernardo, M., Gorochowski, T.E., \& Russo, G. Synchronization and control of complex networks via contraction, adaptation and evolution. IEEE Circuits Syst. Mag. 10, 64-82 (2010).
\bibitem{zhang} Zhang, X., Hu, X., Kurths, J., \& Liu, Z. Explosive synchronization in a general complex network. Phys. Rev. E 88, 010802 (2013).
\bibitem{Yao} Yao, Z., Ma, J., Yao, Y., \& Wang, C. Synchronization realization between two nonlinear circuits via an induction coil coupling. Nonlinear Dyn. 96, 205-217 (2019).
\bibitem{sujith} Sujith, R.I., and Unni, V.R.  Complex system approach to investigate and mitigate thermoacoustic instability in turbulent combustors. Phys. Fluids 32, 061401 (2020).


\bibitem{Pecora1}Pecora, L.M. and Carroll, T.L. Synchronization in chaotic systems. Phys. Rev. Lett. 64, 821 - 824 (1990).
\bibitem{Pecora2}Pecora, L.M. and Carroll, T.L.  Driving systems with chaotic signals. Phys. Rev. A. 44, 2374 - 2383 (1991).
\bibitem{Pecora3} Pecora, L.M., and Carroll, T.L.  Synchronization of chaotic systems. Chaos 25, 097611 (2015).

\bibitem{Boccaletti} Boccaletti, S., Kurths, J., Osipov, G., Valladares, D.L., \& Zhou, C.S.  The synchronization of chaotic systems. Phys. Rep. 366, 1-101 (2002).

\bibitem{Pecora4} Pecora, L., Carroll, T., Johnson, G., Mar, D., \& Fink, K.S.  Synchronization stability in coupled oscillator arrays: Solution for arbitrary configurations. Int. J.  Bifurcat. Chaos 10, 273-290 (2000).

\bibitem{Rulkov:1995} Rulkov, N.F., Sushchik, M.M., Tsimring, L.S. \& Abarbanel, H.D.I. Generalized synchronization of chaos in directionally coupled chaotic systems. Phys. Rev. E51, 980 - 994 (1995). 

\bibitem{kocarev:1996} Kocarev, L. and Parlitz, U.  Generalized synchronization, predictability, and equivalence of unidirectional coupled dynamical systems, Phys. Rev. Lett. 76, 1816 - 1819 (1996).

\bibitem{Ding} Ding, M. and Ott, E.  Enhancing synchronism of chaotic systems, Phys. Rev. E49, R945 (1994).

\bibitem{Kapitaniak} Kapitaniak, T.  Synchronization of chaos using continuous control, Phys. Rev. E 50, 1642 - 1644 (1994).

\bibitem{Gu} Gu, X., Jia, F., Deng, Z., \& Hu, R.  Stochastic response of nonlinear viscoelastic systems with time-delayed feedback control force and bounded noise excitation. Int. J. Struct. Stab. Dyn., 21, 2150181 (2021).

\bibitem{Sayed1}Hamed, Y. S., Albogamy, K. M., \& Sayed, M.  Nonlinear vibrations control of a contact-mode AFM model via a time-delayed positive position feedback. Alex. Eng. J, 60, 963-977 (2021).

\bibitem{Sayed2}Sayed, M., Mousa, A. A., \& Alzaharani, D. Y.  Non-linear time delay saturation controller for reduction of a non-linear vibrating system via 1: 4 internal resonance. J. Vibroengineering, 18, 2515-2536 (2016).

\bibitem{Cantisan} Cantis\'an, J., Coccolo, M., Seoane, J.M., \& Sanju\'an, M.A.F.  Delay-induced resonance in the time-delayed duffing oscillator. Int. J.  Bifurcat. Chaos 30, 2030007 (2020).

\bibitem{Coccolo} Coccolo, M., Cantis\'an, J., Seoane, J.M., Rajasekar, S., \& Sanju\'an, M.A.F. Delay-induced resonance suppresses damping-induced unpredictability. Philos. Trans. Royal Soc. A 379, 20200232 (2021).

\bibitem{Miranda} Gonz\'alez-Miranda, J. M. Synchronization and control of chaos: an introduction for scientists and engineers. World Scientific (2004).





\end{thebibliography}
\end{document}